\documentclass[final,3p]{elsarticle}

\usepackage{amsmath,amsthm,amssymb,mathrsfs}
\usepackage{enumerate}
\usepackage{slashed}
\usepackage{txfonts,dsfont}
\usepackage{aliascnt}
\usepackage[colorlinks, linkcolor=red, citecolor=green,filecolor=magenta,urlcolor=cyan]{hyperref}

\allowdisplaybreaks

\makeatletter
\newcommand{\autorefcheckize}[1]{%
  \expandafter\let\csname @@\string#1\endcsname#1%
  \expandafter\DeclareRobustCommand\csname relax\string#1\endcsname[1]{%
    \csname @@\string#1\endcsname{##1}\wrtusdrf{##1}}%
  \expandafter\let\expandafter#1\csname relax\string#1\endcsname
}
\makeatother

\allowdisplaybreaks
\newcommand{\abs}[1]{\left\lvert#1\right\rvert}
\newcommand{\kh}[1]{\left(#1\right)}%
\newcommand{\zkh}[1]{\left[#1\right]}%
\newcommand{\dkh}[1]{\left\{#1\right\}}%
\newcommand{\jkh}[1]{\left\langle #1\right\rangle}

\newcommand{\Rmn}[1]{\uppercase\expandafter{\romannueral#1}}

\theoremstyle{plain}
\newtheorem{theorem}{Theorem}[section]

\newaliascnt{lem}{theorem}
\newtheorem{lem}[lem]{Lemma}
 \aliascntresetthe{lem}

\newaliascnt{cor}{theorem}
\newtheorem{cor}[cor]{Corollary}
 \aliascntresetthe{cor}

\newaliascnt{prop}{theorem}

 \aliascntresetthe{prop}

\newtheorem*{conj}{Conjecture}

\newtheorem{question}{Question}

\theoremstyle{remark}
\newtheorem{rem}{Remark}[section]

\theoremstyle{definition}

\numberwithin{equation}{section}

\journal{a journal}

\begin{document}

\begin{frontmatter}


\title{A new characterization for Clifford hypersurfaces}

\author[cui]{Qing Cui}
\address[cui]{School of Mathematics,
Southwest Jiaotong University, 611756
Chengdu, Sichuan, China}
\ead{cuiqing@swjtu.edu.cn}

\author[pen]{Carlos Pe\~{n}afiel}
\address[pen]{Instituto de Matem\'{a}tica, Universidade Federal de Rio de Janeiro, Rio de Janeiro, 
22453-900, Brazil}

\ead{penafiel@im.ufrj.br
}

\date{\today ~ [120430]}

\begin{abstract}
 For a closed minimal immersed hypersurface $M$ in $\mathbb S^{n+1}$ with second fundamental form $A$, and each integer $k\ge 2$, define a constant $\sigma_k=\dfrac{\int_M \kh{\abs{A}^2}^k}{\abs{M}}$. We show that $\sigma_k \ge 2^k$ provided $n=2$ and $M$ is not totally geodesic. When $n=4$ and $M$ has two distinct principal curvatures, we show $\sigma_2 \ge 16$. When $n\ge 3$ and $M$ has two distinct principal curvatures, for each integer $k\ge 2$, there exists a positive constant $\delta_k(n)<n$, if $\abs{A}^2\ge \delta_k(n)$, we have $\sigma_k\ge n^k$. All the equality holds iff $M$ is isometric to a Clifford hypersurface. 
\end{abstract}


\end{frontmatter}

\section{Introduction}

Let $M$ be a closed minimal surface in 3-sphere $\mathbb S^3$ with second fundamental form $A$ and Gaussian curvature $K$, then the Gauss equation reads
\begin{align*}
K=1-\dfrac{1}{2} \abs{A}^2.
\end{align*}
Integrating the above equality on $M$ and using Gauss-Bonnet formula, one has
\begin{align}\label{GBdim2}
\sigma:=\dfrac{\int_M \abs{A}^2}{\abs{M}}= 2-\dfrac{8\pi(1-g_M)}{\abs{M}},
\end{align}
where $g_M$ denotes the genus of $M$. If $M$ is not totally geodesic, then by Almgren's uniqueness result \cite{Alm66}, one has $g_M\ge 1$. Consequently, $\sigma\ge 2$. Moreover, if we assume $M$ is embedded, by Lawson's conjecture which is proved by Brendle \cite{Bre13a}, $\sigma=2$ iff $M$ is a Clifford torus. It is natural to ask whether these results still hold for higher dimension. To this end, in 2004, Perdomo \cite{Per04} conjectured that 
\begin{conj}[Perdomo Conjecture]
Let $M$ be a closed embedded non-totally geodesic minimal hypersurface in $\mathbb S^{n+1}$, then $\sigma \ge n$. Moreover, $\sigma =n$ iff $\abs{A}^2\equiv n$ and $M$ is a Clifford torus. 
\end{conj}
Very recently, Ge and Li \cite{GeLi22} gave a partial affirmative answer to Perdomo's conjecture. They showed that $\sigma \ge \delta(n)$ if $M$ satisfies the assumption of Perdomo's conjecture, where $\delta(n)$ is a positive constant depends only on $n$. Note that the assumption ``embedded'' is necessary for Perdomo's conjecture. Actually, Perdomo himself \cite{Per04a} showed that if $M$ is a closed immersed minimal hypersurface in $\mathbb S^{n+1}$ with two distinct principal curvatures, then the inverse inequality holds, i.e., $\sigma \le n$.

In this paper, for a closed immersed minimal hypersurface $M$ in $\mathbb S^{n+1}$ with second fundamental form $A$, we consider the following constants:
\begin{align}\label{sigma_k}
\sigma_k:=\dfrac{\int_M \kh{\abs{A}^2}^k}{\abs{M}}, \quad \text{for  integer}\ \ k\ge 2.
\end{align}
Constants $\sigma_k$ can be used to characterize the Clifford tori. We will summarize our results in the following main theorem:
\begin{theorem}
Let $M$ be a closed immersed minimal non-totally geodesic hypersurface in $\mathbb S^{n+1}$ with $\sigma_k$ defined by \eqref{sigma_k}, then we have
\begin{itemize}
\item $n=2$, then $\sigma_k\ge 2^k$, with equality holds iff $M$ is isometric to $\mathbb S^1\kh{\dfrac{1}{\sqrt{2}}}\times \mathbb S^1\kh{\dfrac{1}{\sqrt{2}}}$;
\item $n=4$ and $M$ has two distinct principal curvatures, then $\sigma_2\ge 16$, with equality holds iff $M$ is isometric to one of the Clifford hypersurfaces $\mathbb S^{l}\kh{\sqrt{\dfrac{l}{4}}}\times \mathbb S^{4-l}\kh{\sqrt{\dfrac{4-l}{4}}}$, for $l=1,2$;
\item $n\ge 3$ and $M$ has two distinct principal curvatures. For each integer $k\ge 2$, there exists constant $\delta_k(n)$ (decreasing on $k$ and $\delta_k(n)\le\delta_1(n)=\dfrac{n(n-2)}{n+2}$), if $\abs{A}^2\ge \delta_k(n)$, then $\sigma_k\ge n^k$. The equality holds iff $M$ is a Clifford hypersurface $\mathbb S^{m}\kh{\sqrt{\dfrac{m}{n}}}\times \mathbb S^{n-m}\kh{\sqrt{\dfrac{n-m}{n}}} (1\le m\le n-1)$.
\end{itemize}
\end{theorem}

As an application, we can give a new characterization for Clifford hypersurface as follows.
\begin{cor}
Let $M$ be a closed immersed minimal hypersurface in $\mathbb S^{n+1}$ with $\sigma_k$ defined by \eqref{sigma_k}, then we have
\begin{itemize}
\item Assume $n=2$,  if for some $k\ge 2$, $\sigma_k= 2^k$, then $M$ is isometric to $\mathbb S^1\kh{\dfrac{1}{\sqrt{2}}}\times \mathbb S^1\kh{\dfrac{1}{\sqrt{2}}}$;
\item Assume $n=4$ and $M$ has two distinct principal curvatures, then $\sigma_2 = 16$ iff $M$ is isometric to one of the Clifford hypersurfaces $\mathbb S^{l}\kh{\sqrt{\dfrac{l}{4}}}\times \mathbb S^{4-l}\kh{\sqrt{\dfrac{4-l}{4}}}$, for $l=1,2$;
\item Assume $n\ge 3$ and $M$ has two distinct principal curvatures. If for some $k\ge 2$, $\abs{A}^2\ge \delta_k(n)$ and $\sigma_k=n^k$, then $M$ is a Clifford hypersurface $\mathbb S^{m}\kh{\sqrt{\dfrac{m}{n}}}\times \mathbb S^{n-m}\kh{\sqrt{\dfrac{n-m}{n}}} (1\le m\le n-1)$.
\end{itemize}
\end{cor}

This paper is organized as follows. In Section 2, we  prove the two special cases of $n=2$ and $n=4$; In Section 3, we prove the general cases. In Section 4, we give some remarks and list several related problems.

\section{$n=2$ and $n=4$ cases}
Gauss-Bonnet-Chern formula will play crucial roles in both of the proof of these two cases. We will state and prove these two cases respectively.
 \begin{theorem}
 Let $M$ be a closed immersed non-totally geodesic minimal surface in $\mathbb S^3$ with $\sigma_k$ defined by \eqref{sigma_k}, then $\sigma_k\ge 2^k$ and the equality holds iff $M$ is isometric to the Clifford torus $\mathbb S^1\kh{\sqrt{\dfrac{1}{2}}}\times \mathbb S^1\kh{\sqrt{\dfrac{1}{2}}}$.
 
 \end{theorem}
\begin{proof}
Denote by $K$ the Gaussian curvature and by $A$ the second fundamental form of $M$, the Gauss equation and Gauss-Bonnet formula yield
\begin{align*}
\sigma_k=\dfrac{\int_M \kh{\abs{A}^2}^k}{\abs{M}}
=\dfrac{\int_M\kh{2-2K}^k}{\abs{M}}=2^k\kh{1-k\dfrac{\int_M K}{\abs{M}}+ \dfrac{\int_M \sum_{i=2}^k C_k^i (-K)^i}{\abs{M}}}=2^k\kh{1+\dfrac{4\pi k(g_M-1)}{\abs{M}}+ \dfrac{\int_M \sum_{i=2}^k C_k^i (-K)^i}{\abs{M}}},
\end{align*}
where $g_M$ is the genus of $M$.
Note that 
\begin{align*}
\sum_{i=2}^k C_k^i (-K)^i=\kh{1-K}^k-1-kK,
\end{align*}
which is a polynomial of $K$. It is easy to see this polynomial take minimum at $K=0$, and consequently, $\sum_{i=2}^k C_k^i (-K)^i\ge 0$. Since $M$ is not totally geodesic, again by Almgren's uniqueness result \cite{Alm66}, $g_M\ge 1$, and therefore $\sigma_k\ge 2^k$. If the equality holds, we have $g_M=1$ and $K\equiv 0$, that is, $M$ is a flat torus, which is isometric to $\mathbb S^1\kh{\sqrt{\dfrac{1}{2}}}\times \mathbb S^1\kh{\sqrt{\dfrac{1}{2}}}$. 
\end{proof}
\begin{rem}
By \eqref{GBdim2}, $M$ has genus one iff $\sigma_1=\sigma=2$. But when $M$ is immersed, there are infinite minimal surfaces in $\mathbb S^3$ with genus one. Therefore, $\sigma=2$ does not characterize the Clifford torus. However, when $k\ge 2$, via the above theorem, $\sigma_k=2^k$ can characterize the Clifford torus.
\end{rem}

When $n=4, k=2$, we have 
\begin{theorem}\label{dim4thm}
Let $M$ be a closed immersed minimal hypersurface in $\mathbb S^5$ with two distinct principal curvatures, then $\sigma_2\ge 16$, i.e., $\int_M \abs{A}^4 \ge 16\abs{M}$. Moreover, the equality holds iff $M$ is isometric to $\mathbb S^{l}\kh{\sqrt{\dfrac{l}{4}}}\times \mathbb S^{4-l}\kh{\sqrt{\dfrac{4-l}{4}}}$, for $l=1,2$.
\end{theorem}
Before prove Theorem \ref{dim4thm}, we need three lemmas.

\begin{lem}\label{lem1}
Let $M$ be a closed minimal hypersurface in $\mathbb S^5$ with second fundamental form $A$, Euler characteristic $\chi(M)$ and principal curvatures $\dkh{\lambda_i}_{i=1}^4$, then we have
\begin{align}\label{GB2nd}
\int_M \kh{\dfrac{3}{2}\abs{A}^4 -3\sum_{i=1}^4 \lambda_i^4-2\abs{A}^2+12}=16\pi^2\chi\kh{M}.
\end{align}
\end{lem}
\begin{proof}
For a closed 4-dimensional manifold, the Gauss-Bonnet-Chern formula (see \cite{Ave63} or \cite{Bes87}) reads
\begin{align}\label{GBC}
\int_M \kh{\frac{s^2}{3} -\abs{\rm Ric}^2+\frac{\abs{W}^2}{2}}=16\pi^2 \chi(M),
\end{align}
where $s, Ric, W$ are the scalar curvautre, Ricci tensor, Weyl curvature tensor of $M$ respectively. Choose a locally orthonormal frame $\dkh{e_i}_{i=1}^n$ on $M$, $\dkh{\omega}_{i=1}^n$ is the dual frame, and assume the second fundamental form $A=h_{ij}\omega^i\otimes \omega^j$, the curvature tensor is $R=R_{ijkl}\omega^i\otimes \omega^j\otimes \omega^k\otimes \omega^l$, then the Gauss equation reads
\begin{align}\label{Gausseq}
R_{ijkl}=\delta_{ik}\delta_{jl} -\delta_{il}\delta_{jk} +h_{ik}h_{jl}-h_{il}h_{jk}.
\end{align}
The following formulas are well-known (see \cite[p.117]{Aub98})
\begin{align*}
s=&\sum_{i,j} R_{ijij},\\
\abs{{\rm Ric}}^2 =&\sum_{i,j} R_{ij}^2=:\sum_{i,j} \kh{\sum_k R_{ikjk}}^2,\\
\abs{W}^2=& \sum_{i,j,k,l} W_{ijkl}^2=\sum_{i,j,k,l} \zkh{R_{ijkl}-\dfrac{1}{2}\kh{R_{ik}\delta_{jl}-R_{il}\delta_{jk}+R_{jl}\delta_{ik}-R_{jk}\delta_{il}}+\dfrac{s}{6}\kh{\delta_{ik}\delta_{jl}-\delta_{il}\delta_{jk}}}^2.
\end{align*}
Note that at a point, $A$ can be diagonalized, i.e., $h_{ij}=\lambda_i \delta_{ij}$, and the minimality of $M$ implies $\sum_i \lambda_i=0$. Then combined the above formulas and Gauss equation \eqref{Gausseq}, a direct computation yields
\begin{align*}
s=&12-\abs{A}^2,\\
\abs{Ric}^2=&36-6\abs{A}^2+\sum_i\lambda_i^4,\\
\abs{W}^2=&\dfrac{7}{3}\abs{A}^4-4\sum_i\lambda_i^4.
\end{align*}
Substituting the above formulas into Gauss-Bonnet-Chern formula \eqref{GBC}, we obtain \eqref{GB2nd}.
\end{proof}

\begin{rem}
Since for $n\ge 4$, locally conformally flat is equivalent to $W\equiv 0$. Therefore, by the proof of the above lemma, a minimal hypersurface in $\mathbb S^5$ is locally conformally flat iff $\abs{A}^4=\dfrac{12}{7}\sum_i \lambda_i^4$. By a direct computation, the trace-free Ricci tensor $\mathring{Ric}=Ric -\dfrac{s}{4}Id$ satisfies
\begin{align*}
\abs{\mathring{Ric}}^2=\sum_i\lambda_i^4 -\dfrac{1}{4}\abs{A}^4.
\end{align*}
Hence, $M$ is Einstein iff $\mathring{Ric}\equiv 0$, i.e., $4\sum_i\lambda_i^4 =\abs{A}^4$.
\end{rem}
The following result is due to \^{O}tsuki.
\begin{lem}[\cite{Ots70}]\label{Otsuki1}
Let $M$ be a closed mininal hypersurface in $\mathbb S^{n+1}$ with two distinct principal curvatures whose multiplicities are both at least two, then $M$ is a Clifford hypersurface $\mathbb S^{m}\kh{\sqrt{\dfrac{m}{n}}}\times \mathbb S^{n-m}\kh{\sqrt{\dfrac{n-m}{n}}} (2\le m\le n-2)$.
\end{lem}
We also need a result due to Perdomo.
\begin{lem}[\cite{Per04a}]\label{Perlem}
Let $M$ be a closed minimal hypersurface in $\mathbb S^{n+1}$ with second fundamental form $A$ and two distinct principal curvatures, then we have $\sigma \le n$, i.e.,
$
\int_M \abs{A}^2 \le n\abs{M}.
$
Moreover, the equality holds iff $M$ is a Clifford hypersurface.
\end{lem}

Now we can give the proof of  Theorem \ref{dim4thm}.

\begin{proof}[Proof of Theorem \ref{dim4thm}]
Assume $M$ has two distinct principal curvatures $\lambda, \mu$. By Lemma \ref{Otsuki1}, if the multiplicities of $\lambda, \mu$ are both two, then $M$ is isometric to $\mathbb S^2\kh{\sqrt{\dfrac{1}{2}}}\times \mathbb S^2\kh{\sqrt{\dfrac{1}{2}}}$. Consequently, $\abs{A}^2\equiv 4$, and the conclusion follows.

Therefore, we next assume $\mu<0$ has multiplicity one and $\lambda>0$ has multiplicity 3. In this case, we have $\mu=-3\lambda$ and 
\begin{align*}
\abs{A}^2 =& (-3\lambda)^2+3\lambda^2=12\lambda^2,\\
\sum_i\lambda_i^4 =& (-3\lambda)^4+3\lambda^4 = 84\lambda^4 =\dfrac{7}{12}\abs{A}^4.
\end{align*}
Hence, by Lemma \ref{lem1} and Lemma \ref{Perlem}, we have 
\begin{align*}
16\pi^2\chi\kh{M}=&\int_M \kh{\dfrac{3}{2}\abs{A}^4 -3\sum_{i=1}^4 \lambda_i^4-2\abs{A}^2+12}\\
=&\int_M \kh{-\dfrac{1}{4}\abs{A}^4 -2\abs{A}^2+12}\\
\ge &\int_M \kh{-\dfrac{1}{4}\abs{A}^4 +4}.
\end{align*}
On the other hand, since in our case $\mu$ is the eigenvalue of $A$ has multiplicity one everywhere, the unit eigenvectors of $\mu$ form a nowhere vanishing vector field on $M$, which implies $\chi(M)=0$ by the Poincar\'{e}-Hopf theorem. Therefore, we have $\int_M \abs{A}^4\ge 16\abs{M}$. If the equality holds, check about the inequality in the proof, we have $\int_M \abs{A}^2 = 4\abs{M}$ which implies, by Lemma \ref{Perlem}, $M$ is the Clifford hypersurface $\mathbb S^1\kh{\dfrac{1}{2}}\times \mathbb S^3\kh{\dfrac{\sqrt{3}}{2}}$. 
\end{proof}

\section{General case}
In this section, we will prove the following theorem.
\begin{theorem}\label{generalthm}
Let $M$ be a closed minimal hypersurface in $\mathbb S^{n+1} (n\ge 3)$ with second fundamental form $A$ and  two distinct principal curvatures. For each integer $k\ge 2$, there exists constant $\delta_k(n)$ (decreasing on $k$ and $\delta_k(n)\le\delta_1(n)=\dfrac{n(n-2)}{n+2}$), if $\abs{A}^2\ge \delta_k(n)$, then $\sigma_k\ge n^k$. The equality holds iff $M$ is a Clifford hypersurface.
\end{theorem}

\begin{proof}

By Lemma \ref{Otsuki1}, in what follows, we only focus on the case that $M$ has two distinct principal curvatures $\mu<0$ with multiplicity one, and $\lambda>0$ with multiplicity $n-1$. In this case, \^{O}tsuki \cite{Ots70} showed that the distribution
\begin{align*}
\Gamma=\dkh{v\in TM \ \vert  A(v)=\lambda v}
\end{align*}
 is completely integrable and $\lambda$ is constant on each integral submanifold of the distribution.
Therefore, if we assume $\dkh{e_i}_{i=1}^n$ be a local orthonormal frame such that $e_n$ is an eigenvector of $A$ w.r.t. $\mu$, i.e., $A(e_n)=\mu e_n$. Then we have 
\begin{align*}
e_i(\mu)=0, \ \ e_i(\lambda)=0, \ \ i=1,\cdots, n-1.
\end{align*}
This implies $\lambda $ and $\mu$ are locally only depend on $e_n$.
For the sake of simplicity, we write $e_n(\lambda)=\dot{\lambda}, e_n(e_n(\lambda))=\ddot{\lambda}$, and so on.  \^{O}tsuki also showed 
$\lambda$ satisfies a second order ODE:
\begin{align}\label{OtsODE}
\ddot{\lambda}-\dfrac{n+1}{n\lambda}\kh{\dot{\lambda}}^2+n\lambda\zkh{(n-1)\lambda^2-1}=0.
\end{align}
For each $i=1,\cdots, n-1$, Perdomo \cite{Per04a} showed that
\begin{align}\label{Per1}
\jkh{\nabla_{e_i}e_n, e_i}=-\dfrac{\dot{\lambda}}{n\lambda}.
\end{align}
Let $f$ be a smooth function of $\lambda$, denote by 
\begin{align*}
f^\prime=\dfrac{df}{d\lambda}, \quad f^{\prime\prime}=\dfrac{d^2 f}{d\lambda^2}, \quad \dot{f}:=e_n(f(\lambda))=f^\prime\dot{\lambda},\quad \ddot{f}:=e_n\kh{e_n(f(\lambda))}=f^{\prime\prime}\kh{\dot{\lambda}}^2 +f^\prime \ddot{\lambda}.
\end{align*}
Therefore, by \eqref{OtsODE} and \eqref{Per1}, we have
\begin{align*}
\Delta f=&\sum_i\jkh{\nabla_{e_i}\nabla f, e_i}\\
=&\sum_i\jkh{\nabla_{e_i}\dot{f}e_n, e_i}\\
=&\ddot{f} +\dot{f}\dfrac{(n-1)\dot{\lambda}}{-n\lambda}\\
=&f^{\prime\prime}\kh{\dot{\lambda}}^2+f^\prime\ddot{\lambda}+f^\prime\dot{\lambda}\dfrac{(n-1)\dot{\lambda}}{-n\lambda}\\
=&f^{\prime\prime}\kh{\dot{\lambda}}^2+f^\prime\kh{\dfrac{n+1}{n\lambda}\kh{\dot{\lambda}}^2-n\lambda\kh{(n-1)\lambda^2-1}}+f^\prime\dot{\lambda}\dfrac{(n-1)\dot{\lambda}}{-n\lambda}\\
=&\kh{f^{\prime\prime}+\dfrac{2}{n\lambda}f^\prime}\kh{\dot{\lambda}}^2-n\lambda\kh{(n-1)\lambda^2-1}f^\prime.
\end{align*}
Integrating the above equality on $M$, we have,
\begin{align}\label{keyeq}
\int_M n\lambda\kh{(n-1)\lambda^2-1}f^\prime =\int_M\kh{f^{\prime\prime}+\dfrac{2}{n\lambda}f^\prime}\kh{\dot{\lambda}}^2.
\end{align}
For each integer $k\ge 2$, taking 
\begin{align*}
f(\lambda)=n^{k-1}\kh{\ln \lambda+\sum_{i=1}^{k-1}\dfrac{(n-1)^i}{2i}\lambda^{2i}}
\end{align*}
in equality \eqref{keyeq}. We obtain
\begin{align*}
f^\prime=&n^{k-1}\kh{\dfrac{1}{\lambda}+\sum_{i=1}^{k-1}\kh{n-1}^i\lambda^{2i-1}},\\
f^{\prime\prime}=&n^{k-1}\kh{-\dfrac{1}{\lambda^2}+\sum_{i=1}^{k-1}\dfrac{\kh{n-1}^i}{2i-1}\lambda^{2(i-1)}}.
\end{align*}
Consequently, \eqref{keyeq} becomes,
\begin{align*}
\int_M\kh{\abs{A}^{2k}-n^k}=&\int_M n\lambda\kh{(n-1)\lambda^2-1}n^{k-1}\kh{\dfrac{1}{\lambda}+\sum_{i=1}^{k-1}\kh{n-1}^i\lambda^{2i-1}}\\
=&\int_M n^{k-1}\zkh{\sum_{i=1}^{k-1}(n-1)^i\kh{\dfrac{1}{2i-1}+\dfrac{2}{n}}\lambda^{2(i-1)}-\dfrac{n-2}{n\lambda^2}}\kh{\dot{\lambda}}^2.
\end{align*}
We assume the positive root of 
\begin{align*}
\sum_{i=1}^{k-1}(n-1)^i\kh{\dfrac{1}{2i-1}+\dfrac{2}{n}}\lambda^{2(i-1)}-\dfrac{n-2}{n\lambda^2}
\end{align*}
is $n(n-1)\kh{\delta_k(n)}^2$ for $\delta_k(n)>0$ only depends on $k$ and $n$. It is easy to see $\delta_k(n)$ is decreasing on $k$ and $\delta_k(n)\le \delta_1(n)=\dfrac{n(n-2)}{n+2}<n$.
Then if $\abs{A}^2\ge \delta_k(n)$, we have 
\begin{align*}
\int_M\kh{\abs{A}^{2k} -n^k} \ge 0,
\end{align*}
that is $\sigma_k\ge n^k$. If the equality holds, we have 
\begin{align*}
\abs{A}^2\equiv \delta_k(n), \quad \text{or }\quad \kh{\dot{\lambda}}^2\equiv 0.
\end{align*}
Since $0<\delta_k(n)<n$, $\abs{A}^2\equiv \delta_k(n)$ can not hold (see \cite{Sim68}). Therefore, we have 
\begin{align*}
\abs{\nabla A}^2=n(n-1) \kh{\dot{\lambda}}^2\equiv 0,
\end{align*}
which implies (see \cite{ChedoCKob68}) $M$ is a Clifford torus.
\end{proof}

\section{Remarks and questions}
\begin{question}
Let $M$ be a closed minimal surface in $\mathbb S^3$, with second fundamental form $A$, for each positive integer $k$, write $\sigma_k =\dfrac{\int_M\kh{\abs{A}^2}^k}{\abs{M}}$. Apparently, $M$ is totally geodesic iff $\sigma_k=0$. If $\sigma_k>0$, then $\sigma_k\ge 2^k$. A natural question is: what is the next value of $\sigma_k$?
\end{question}
When $k=1$ and $M$ is embedded, 
if $\sigma_1>2$, then uniqueness of minimal embedded torus (see \cite{Bre13a}) and \eqref{GBdim2} imply that the genus $g_M\ge 2$. Therefore, 
the area estimate of minimal surface in $\mathbb S^3$ (see Choi and Wang's result \cite{ChoWan83})  and \eqref{GBdim2} yield,
\begin{align*}
\sigma_1=  2+\dfrac{8\pi(g_M-1)}{\abs{M}}\ge 2+\dfrac{8\pi(g_M-1)}{8\pi (g_M+1)}\ge \dfrac{7}{3}.
\end{align*}
Since the area estimate of Choi and Wang based on the estimate of the first positive eigenvalue of Laplacian, the first inequality in the above formula is actually strict (see \cite[Theorem 5.1]{Bre13b}). Therefore, we have if $\sigma_1>2$ then $\sigma_1>\dfrac{7}{3}$.

For $n=2, k\ge 2$, the inequality $\sigma_k\ge n^k$ always hold. It is natural to ask:
\begin{question}
Let $M$ be a closed non-totally geodesic minimal hypersurface in $\mathbb S^{n+1}$. Is the inequality $\sigma_k\ge n^k$  hold? That is, is the inequality
\begin{align*}
\int_M \abs{A}^{2k} \ge n^k\abs{M}, \ \ \text{for} \ \ k\ge 2
\end{align*}
hold?
\end{question}
Note that, for $k= 1$ the inequality does not hold for general immersed case \cite{Per04a}.  But for $k\ge 2$, our main theorem gives a partial positive answer to the above question. In particular, when $k=2$,  the Simons' identity (see \cite{Sim68}) gives,
\begin{align*}
\dfrac{1}{2}\Delta \abs{A}^2=\abs{\nabla A}^2+(n-\abs{A}^2)\abs{A}^2.
\end{align*}
Integrating the above equality on $M$, we get,
\begin{align*}
\int_M \abs{A}^4=\int_M\abs{\nabla A}^2 +n\int_M \abs{A}^2 \ge n\int_M \abs{A}^2.
\end{align*}
Moreover, if we assume $M$ is embedded, and if the Perdomo's conjecture hold, we have $\int_M\abs{A}^4\ge n^2\abs{M}$, i.e., $\sigma_2\ge n^2$.

\providecommand{\bysame}{\leavevmode\hbox to3em{\hrulefill}\thinspace}
\providecommand{\MR}{\relax\ifhmode\unskip\space\fi MR }

\end{document}